\documentclass[11pt,a4paper]{article}
\usepackage{graphicx}
\usepackage{amsfonts}
\usepackage{latexsym}
\usepackage{epsfig}

\newcommand{\beql}[1]{\begin{equation}\label{#1}}
\newcommand{\eeq}{\end{equation}}
\newcommand{\comment}[1]{}

\newcommand{\eqref}[1]{{\rm (\ref{#1})}}

\newcommand{\Abs}[1]{{\left|{#1}\right|}}

\newcommand{\Norm}[1]{{\left\|{#1}\right\|}}

\newcommand{\Qed}{\ \\\mbox{$\Box$}}

\newcommand{\Set}[1]{{\left\{{#1}\right\}}}

\newcommand{\RR}{{\mathbb R}} 

\newcommand{\ZZ}{{\mathbb Z}}

\newcommand{\inner}[2]{{\langle #1, #2 \rangle}}

\newcommand{\dens}{{\rm dens\,}}

\newcommand{\supp}{{\rm supp\,}}
\newcommand{\dist}{{\rm dist\,}}

\newcommand{\ft}[1]{\widehat{#1}}

\newcounter{open}
\newcounter{dfn}
\def\thedfn{\arabic{dfn}}
\newenvironment{dfn}{
  \sf
  \vskip 0.10in
  \refstepcounter{dfn}
  \noindent{\bf Definition \thedfn \ }
}{\vskip 0.10in}

\newcounter{obs}
\def\theobs{\arabic{obs}}

\newcounter{thm}
\newcounter{othm}
\def\theothm{\Alph{othm}}
\newenvironment{othm}{
  \sf
  \vskip 0.10in
  \refstepcounter{othm}
  \noindent{\bf Theorem\ \theothm}
}{\vskip 0.10in}

\newcounter{mysec}
\def\themysec{\arabic{mysec}}
\newcommand{\mysection}[1]{
  \vskip 0.25in
  \refstepcounter{mysec}\noindent{\large\bf \S\themysec.\ {#1}}\par\ \par
  \nopagebreak
  \addcontentsline{toc}{section}{{\bf \themysec.}\ {#1}}
}

\newcounter{mysubsec}[mysec]
\newtheorem{theorem}{Theorem}
\newtheorem{corollary}{Corollary}
\newtheorem{lemma}{Lemma}

\newtheorem{conjecture}{Conjecture}

\begin{document}

\title{On a problem of Tur\'an about positive definite functions}
\author{Mihail N. Kolountzakis \and Szil\' ard Gy. R\' ev\' esz
\thanks {The second author was supported in part by the 
Hungarian National Foundation for Scientific Research, 
Grant \# T034531 and T 032872.}
}
\date{April 2002}
\maketitle

\begin{abstract}
We study the following question posed by Tur\' an. Suppose
$\Omega$ is a convex body in Euclidean space $\RR^d$ which is symmetric
with respect to the origin. Of all positive definite functions supported
in $\Omega$, and with value $1$ at the origin, which one has the largest integral?
It is probably the case that the extremal function is the indicator of the
half-body convolved with itself and properly scaled, but this has been proved
only for a small class of domains so far.
We add to this class of known
{\em Tur\' an domains} the class of all
spectral convex domains. These are all convex domains which have
an orthogonal basis of exponentials $e_\lambda(x) = \exp 2\pi i\inner{\lambda}{x}$,
$\lambda \in \RR^d$.
As a corollary we obtain that all convex domains which tile space by translation
are Tur\' an domains.

We also give a new proof that the Euclidean ball is a Tur\' an domain.
\end{abstract}



%%%%%%%%%%%%%%%%%%%%%%%%%%%%%%%%%%%%%%%%%%%%%
% Description of the TURAN PROBLEM
%%%%%%%%%%%%%%%%%%%%%%%%%%%%%%%%%%%%%%%%%%%%%
\mysection{A problem of Tur\'an}

\noindent
The following question is attributed to Tur\' an: 
\begin{quotation}
({\bf Tur\' an}) Let $\Omega$ be a convex domain in $\RR^d$ which is
symmetric with respect to $0$.
What is the maximum of $\int f$ for all $f$ supported in
$\Omega$ which are positive definite (their Fourier Transform is nonnegative)
and have $f(0)=1$?
\end{quotation}
See \cite{hexagon,ball,stechkin} for the history of the problem.

\begin{dfn} (Tur\' an domains)\\
A symmetric convex domain $\Omega$  is called a {\em  Tur\' an domain} if
the maximum of $\int f$ over all positive definite functions
supported in $\Omega$ and with $f(0)=1$, is achieved by the function
\beql{half}
f = \frac {2^d}{|\Omega|} \chi_{{1\over2}\Omega}*\chi_{{1\over2}\Omega}.
\eeq
Otherwise it is called a {\em non-Tur\' an domain}.
\end{dfn}
The function $f$ in the above definition is clearly a legitimate function,
as its Fourier Transform (FT) is
$$
\widehat{f} = \frac{2^d}{|\Omega|} \Abs{\ft{\chi_{{1\over2}\Omega}}}^2
$$
and $f(0)=1$.
The integral of $f$ in this case is $2^{-d}|\Omega|$.

Non-Tur\' an domains are not known to exist.

We cite the recent papers by Arestov and Berdysheva
\cite{hexagon} and by Gorbachev \cite{ball}.
In \cite{hexagon} it is proved that the regular
hexagon in the plane is a Tur\' an domain
and in \cite{ball} that the Euclidean ball
in $d$-dimensional space is also
one. We find it very likely that there are no non-Tur\' an domains.

\noindent
{\em Remark.} This problem is only interesting for {\em symmetric} domains as
the support of any positive definite function is always a
symmetric set. The convexity assumption on the domain also cannot be removed 
if the function $f$ in \eqref{half} is to be supported in $\Omega$.

\mysection{A problem of Fuglede}

Our interest in this problem of Tur\' an arose in connection with a
conjecture of Fuglede which we now describe.

%%%%%%%%%%%%%%%%%%%%%%%%%%%%%%%%%%%%%%%%%%%%%
% Define TILES
%%%%%%%%%%%%%%%%%%%%%%%%%%%%%%%%%%%%%%%%%%%%%
\begin{dfn} (Translational tiles)\\ 
A measurable set $\Omega\subseteq\RR^d$ is a {\em translational tile} 
if there exists a set $\Lambda\subseteq\RR^d$ such that almost all (Lebesgue) 
points in $\RR^d$ belong to exactly one of the translates 
$$ 
\Omega + \lambda,\ \ \ \lambda\in\Lambda. 
$$ 
We denote this condition by $\Omega + \Lambda = \RR^d$. 

If $f\in L^1(\RR^d)$ is nonnegative we say that $f$ tiles with
$\Lambda$ at level $\ell$ if
$$
\sum_{\lambda\in\Lambda} f(x-\lambda) = \ell,\ \ \mbox{\rm a.e.\ $x$}.
$$
We denote this latter condition by $f+\Lambda = \ell\RR^d$.
\end{dfn} 

%%%%%%%%%%%%%%%%%%%%%%%%%%%%%%%%%%%%%%%%%%%%%
% UNIFORM DENSITY
%%%%%%%%%%%%%%%%%%%%%%%%%%%%%%%%%%%%%%%%%%%%%
In any tiling the translation set has some properties of
density, which hold uniformly in space.
\begin{dfn} \label{dfn:density} (Uniform density)\\
A multiset $\Lambda\subseteq\RR^d$ has (uniform) density
$\rho$ if
$$
\lim_{R\to\infty} {\#(\Lambda \cap B_R(x)) \over \Abs{B_R(x)}} \to
\rho
$$
uniformly in $x\in\RR^d$. We write $\rho = \dens \Lambda$.\\
We say that $\Lambda$ has (uniformly) bounded density if the fraction
above is bounded by a constant $\rho$ uniformly for $x\in\RR$ and
$R>1$.
We say then that $\Lambda$ has density (uniformly) bounded by $\rho$.
\end{dfn}

\noindent
{\em Remark.}
It is not hard to prove (see for example \cite{line}, Lemma 2.3, where
it is proved in dimension one -- the proof extends verbatim to higher 
dimension)
that in any tiling $f+\Lambda = \ell\RR^d$ the set
$\Lambda$ has density $\ell / \int f$.

%%%%%%%%%%%%%%%%%%%%%%%%%%%%%%%%%%%%%%%%%%%%%
% The L^2 space and sets of exponentials
%%%%%%%%%%%%%%%%%%%%%%%%%%%%%%%%%%%%%%%%%%%%%
Let $\Omega$ be a measurable subset of $\RR^d$ 
and $\Lambda$ be a discrete subset of $\RR^d$.
We write
\begin{eqnarray*}
e_\lambda(x) &=& \exp{2\pi i \inner{\lambda}{x}},\ \ \ (x\in\RR^d),\\
E_\Lambda &=& \Set{e_\lambda:\ \lambda\in\Lambda} \subset L^2(\Omega).
\end{eqnarray*}
The inner product and norm on $L^2(\Omega)$ are
$$
\inner{f}{g}_\Omega = \int_\Omega f \overline{g},
\ \mbox{ and }\
\Norm{f}_\Omega^2 = \int_\Omega \Abs{f}^2.
$$

%%%%%%%%%%%%%%%%%%%%%%%%%%%%%%%%%%%%%%%%%%%%%
% SPECTRAL PAIRS
%%%%%%%%%%%%%%%%%%%%%%%%%%%%%%%%%%%%%%%%%%%%%
\begin{dfn} (Spectrum of a domain, spectral pairs)\\
The pair $(\Omega, \Lambda)$ is called a {\em spectral pair}
if $E_\Lambda$ is an orthogonal basis for $L^2(\Omega)$.
A set $\Omega$ will be called {\em spectral} if there is
$\Lambda\subset\RR^d$ such that
$(\Omega, \Lambda)$ is a spectral pair.
The set $\Lambda$ is then called a {\em spectrum} of $\Omega$.
\end{dfn}

\noindent
{\bf Example:} If $Q_d = (-1/2, 1/2)^d$ is the cube of
unit volume in $\RR^d$ then
$(Q_d, \ZZ^d)$ is a spectral pair, as is well known by the 
ordinary $L^2$ theory of multiple Fourier series.

%%%%%%%%%%%%%%%%%%%%%%%%%%%%%%%%%%%%%%%%%%%%%
% The FUGLEDE CONJCETURE
%%%%%%%%%%%%%%%%%%%%%%%%%%%%%%%%%%%%%%%%%%%%%
The following conjecture is still unresolved, in all dimensions and
both directions.

\begin{conjecture}\label{conj:fuglede}
{\rm ({\bf Fuglede \cite{conjecture}})}
Let $\Omega \subset \RR^d$
be a bounded open set.  Then $\Omega$ is spectral
if and only if
there exists $L \subset \RR^d$ such that $\Omega + L = \RR^d$ is a
tiling.
\end{conjecture}
Note that in Fuglede's Conjecture
%%% changed "tiling set" to "translation set"
no relation is claimed between the translation 
set $L$ and the spectrum $\Lambda$.
This Conjecture is still open for all dimensions.
%%% The Conjecture is wide open in all dimensions! The conjecture
%%% does not concern only convex domains. It is only in the convex
%%% case that it has been solved in dimension 2.
However (see \cite{conjecture,nonsym}), the lattice case
of this conjecture is easy to show. 
In the following result the dual lattice $\Lambda^*$ of 
a lattice $\Lambda$ is defined as usual by 
$$ 
\Lambda^* = \Set{x\in\RR^d:\ \forall \lambda\in\Lambda\ \  
\inner{x}{\lambda} \in \ZZ}.
$$ 

%%%%%%%%%%%%%%%%%%%%%%%%%%%%%%%%%%%%%%%%%%%%%
% LATTICE FUGLEDE Theorem
%%%%%%%%%%%%%%%%%%%%%%%%%%%%%%%%%%%%%%%%%%%%%
\begin{othm}\label{th:lattice-fuglede}
{\rm ({\bf Fuglede \cite{conjecture}})}
The bounded, open domain $\Omega$ admits translational tilings 
by a lattice $\Lambda$
if and only if $E_{\Lambda^*}$ is an orthogonal basis for $L^2(\Omega)$.
\end{othm}

\mysection{Connecting the problems of Tur\'an and Fuglede}

%%%%%%%%%%%%%%%%%%%%%%%%%%%%%%%%%%%%%%%%%%%%%
% MAIN THEOREM
%%%%%%%%%%%%%%%%%%%%%%%%%%%%%%%%%%%%%%%%%%%%%
Our main result is the following.
\begin{theorem}\label{th:main}
Let $\Omega \subseteq \RR^d$ be a symmetric convex domain. 
If $\Omega$ is spectral, then it has to be a Tur\'an domain as well.
\end{theorem}

It follows that convex spectral domains are Tur\'an domains. Indeed, 
convex spectral domains are necessarily symmetric according to the result 
in \cite{nonsym}, and so the above Theorem \ref{th:main} applies.

%%%%%%%%%%%%%%%%%%%%%%%%%%%%%%%%%%%%%%%%%%%%%
% MAIN COROLLARY
%%%%%%%%%%%%%%%%%%%%%%%%%%%%%%%%%%%%%%%%%%%%%
\begin{corollary}\label{cor:main}
Suppose the symmetric convex domain $\Omega\subseteq\RR^d$ is
a translational tile. Then it is a Tur\' an domain.
\end{corollary}
\noindent{\em Proof of Corollary \ref{cor:main}.}
We start with the following result which claims that every convex tile
is also a lattice tile.
%%%%%%%%%%%%%%%%%%%%%%%%%%%%%%%%%%%%%%%%%%%%%
% VENKOV - McMULLEN
%%%%%%%%%%%%%%%%%%%%%%%%%%%%%%%%%%%%%%%%%%%%%
\begin{othm}\label{th:venkov}
{\rm({\bf Venkov \cite{venkov} and McMullen \cite{mcmullen}})}
Suppose that a convex body $K$ tiles space by translation.
Then it is necessarily a symmetric polytope and there is a lattice $L$ such that  
$$ 
K+L = \RR^d. 
$$
\end{othm}
A complete characterization of the tiling polytopes is also among
the conclusions of the Venkov-McMullen Theorem but we do not need it here and
choose not to give the full statement as it would require some more 
definitions.

So, if a 
convex domain is a tile, it is also a lattice tile, hence spectral
by Theorem \ref{th:lattice-fuglede}, and as 
such it is Tur\'an,
by Theorem \ref{th:main}.
\Qed

\noindent
{\em Remark.}
If one wants to avoid using the Venkov-McMullen theorem in the
proof of Corollary \ref{cor:main} one should enhance the assumption
of Corollary \ref{cor:main} to state that $\Omega$ is a lattice tile.

The result of \cite{hexagon} about the hexagon being a Tur\' an domain
is thus a special case of our Corollary \ref{cor:main}, but not the result
in \cite{ball} about the ball being Tur\' an. The ball \cite{ball-distances},
and essentially every smooth convex body \cite{smooth}, is known not to
be spectral, in accordance with Conjecture \ref{conj:fuglede}.

Fuglede's Conjecture for convex domains
is still open except for dimension $d=2$, in which case
it was answered in the affirmative recently, see \cite{two-dim}. 
Thus our Theorem \ref{th:main}  conceivably (though not
very likely) applies to a wider class of convex domains than just convex tiles,
dealt with in Corollary \ref{cor:main}.

%%%%%%%%%%%%%%%%%%%%%%%%%%%%%%%%%%%%%%%%%%%%%
% PROOF of MAIN THEOREM
%%%%%%%%%%%%%%%%%%%%%%%%%%%%%%%%%%%%%%%%%%%%%
\noindent
{\em Proof of Theorem \ref{th:main}.}
The proof of our main theorem relies on Fourier theoretic 
characterizations of translational tiling \cite{nonsym}.
Without loss of generality let us assume from now on that 
$\Omega$ has volume $1$.

Let $\Omega$ have spectrum $\Lambda \subseteq \RR^d$.
This is equivalent to the following (see \cite{nonsym})
\beql{tiling}
\sum_{\lambda\in\Lambda} \Abs{\ft{\chi_\Omega}}^2(x-\lambda) = 1,
\ \ \ \mbox{for a.e.\ $x\in\RR^d$}.
\eeq
That is, $\Abs{\ft{\chi_\Omega}}^2$ {\em tiles} $\RR^d$
with translation set $\Lambda$ at level 1, i.e.\ 
$\Abs{\ft{\chi_\Omega}}^2 + \Lambda = \RR^d$.
It follows (see the Remark after Definition \ref{dfn:density}) 
that $\dens \Lambda = 1/|\Omega| = 1$.

For any given $\Lambda \subset \RR^d$ with bounded density (see Definition
\ref{dfn:density}) we denote by $\delta_\Lambda$ the (infinite) measure
$\sum_{\lambda\in\Lambda}\delta_\lambda$. 
%%% The mass grows polnomially, not linearly
This is a tempered distribution, as the total mass in
a ball of radius $R$ grows polynomially with $R$, and therefore
we can speak of its FT.

%%%%%%%%%%%%%%%%%%%%%%%%%%%%%%%%%%%%%%%%%%%%%
% TILING implies SUPPORT CONDITION
%%%%%%%%%%%%%%%%%%%%%%%%%%%%%%%%%%%%%%%%%%%%%
We shall use the following result from \cite{nonsym}.
\begin{lemma}\label{th:tiling-implies-supp}
{\rm({\bf Kolountzakis \cite{nonsym}})}
Suppose that $f\ge 0$ is not identically $0$, that $f \in L^1(\RR^d)$,
$\widehat{f}\ge 0$ has compact support and $\Lambda\subset\RR^d$.
If $f+\Lambda$ is a tiling then
\beql{sp-cond-1}
\supp \ft{\delta_\Lambda} \subseteq \Set{x\in\RR^d:\ \ft{f}(x) = 0} 
\cup \Set{0}.
\eeq
\end{lemma}

When applied to our case, $f = \Abs{\ft{\chi_\Omega}}^2$ and
$\ft{f} = \chi_\Omega * \chi_\Omega$, and it follows that
\beql{support}
\supp \ft{\delta_\Lambda} \subseteq \Set{0} \cup (\Omega-\Omega)^c =
 \Set{0} \cup (2\Omega)^c.
\eeq

%%%%%%%%%%%%%%%%%%%%%%%%%%%%%%%%%%%%%%%%%%%%%
% SUPPORT CONDITION + epsilon implies TILING
%%%%%%%%%%%%%%%%%%%%%%%%%%%%%%%%%%%%%%%%%%%%%
The necessary support condition \eqref{sp-cond-1} in Lemma
\ref{th:tiling-implies-supp} cannot by itself guarantee that $f$ tiles with
$\Lambda$. The reason is that a tempered distribution, such as
$\ft{\delta_\Lambda}$, which is supported in the zero-set of a function, is
not necessarily killed when multiplied by that function. One has to know some
extra information about the order of the distribution (``what order derivatives
it involves'') versus the degree of vanishing of the function on the support of
the distribution\footnote{
An important special case is when
one knows the distribution to be a measure, as is the case when
$\Lambda$ is either a lattice or fully periodic.
In that case any vanishing of the function will do and
the implication in Lemma \ref{th:tiling-implies-supp}
can essentially be reversed.
}.

In the following partial converse to Lemma
\ref{th:tiling-implies-supp} (see \cite{nonsym}), this problem
is solved as the separation of the supports guarantees infinite order of
vanishing of $f$.

\begin{lemma}\label{th:disjoint-supports} {\rm({\bf Kolountzakis
\cite{nonsym}})} Suppose that $g\in L^1(\RR^d)$,
and that $\Lambda \subset \RR^d$ has uniformly bounded density.
Suppose also that $O \subset \RR^d$ is open, that\footnote{
Condition \eqref{non-zero-integral} was mistakenly ommitted from
\cite{nonsym}. This does not affect the validity of what's proved in
\cite{nonsym} as the condition \eqref{non-zero-integral}
is easily seen to hold in the specific application.}
%%%%%%%%%%%%%%%%%%%%%%%%%%%%%%%%%%%%%%%%%%%%%%%%%%%%%%%%%%%%%%%%%%%%%%%%%%%%%
\beql{non-zero-integral}
\ft{g}(0)=\int g \neq 0,
\eeq
and that for some $\delta>0$
\beql{disjoint-supports}
\supp \ft{\delta_\Lambda} \setminus \Set{0} \ \subseteq\ O
\ \mbox{ and }\ O+B_\delta(0) \subseteq\  \Set{\ft{g} = 0}.
\eeq
Then $g+\Lambda$ is a tiling at level
$\ft{g}(0)\cdot\ft{\delta_\Lambda}(\Set{0})$.
\end{lemma}

The conlusion of Lemma \ref{th:disjoint-supports}
demands some explanation.
Conditions \eqref{non-zero-integral} and
\eqref{disjoint-supports} 
imply that in a neighborhood of $0$ the tempered distribution
$\ft{\delta_\Lambda}$ is supported at $0$ only.
That's because $\ft{g}$ is continuous and, since $\ft{g}(0)\neq 0$,
it does not vanish in some neighborhood of $0$.
It then follows that, near $0$, $\ft{\delta_\Lambda}$
is not only a tempered distribution but a measure,
that is, it is just a point mass at $0$ 
(see \cite{line}, Theorem 5.1, Step 1, for the
proof in dimension $1$, which works in any dimension).
For this reason it makes sense to write $\ft{\delta_\Lambda}(\Set{0})$
for that point mass.

From Lemma \ref{th:at-zero} below, it follows that the value of this
constant is precisely the density of $\Lambda$, if such a density exists.

%%%%%%%%%%%%%%%%%%%%%%%%%%%%%%%%%%%%%%%%%%%%%
% VALUE AT ZERO OF FT OF POINT MASSES
%%%%%%%%%%%%%%%%%%%%%%%%%%%%%%%%%%%%%%%%%%%%%
\begin{lemma}\label{th:at-zero}
{\rm ({\bf Kolountzakis \cite{polygons}})}\\
Suppose that $\Lambda \in \RR^d$ is a multiset with density $\rho$,
$\delta_\Lambda = \sum_{\lambda\in\Lambda} \delta_\lambda$,
and that $\ft{\delta_\Lambda}$ is a measure in a neighborhood of
$0$.
Then $\ft{\delta_\Lambda}(\Set{0}) = \rho$.
\end{lemma}

%%%%%%%%%%%%%%%%%%%%%%%%%%%%%%%%%%%%%%%%%%%%%%%%%%%%%%%%%%%
%%%%%      CONCLUSION OF THE PROOF OF THEOREM 2    %%%%%%%%
%%%%%%%%%%%%%%%%%%%%%%%%%%%%%%%%%%%%%%%%%%%%%%%%%%%%%%%%%%%
\noindent
{\em Conclusion of Proof of Theorem \ref{th:main}.}
If $\Omega$ is of non-Tur\'an type, then
there exists a positive definite function $F$ supported 
in $\Omega$ with $F(0) = 1$ and $\int F > 2^{-d}$.

Now define
$$
\ft{G}(x) = F\left({1+\epsilon\over 2} x\right),
$$
where $\epsilon>0$ is to be taken so small that we have $\int \ft{G} > 1$.
It follows that $\supp \ft{G} \subseteq (1-\theta)2\Omega$,
for some $\theta>0$.
The function $\ft{G}$ is also positive definite.

Because of \eqref{support} we can now write
$$
\supp\ft{\delta_\Lambda} \subseteq \Set{0} \cup \Set{\ft{G}=0}.
$$

\begin{figure}\label{fig:domains}
\begin{center}\input{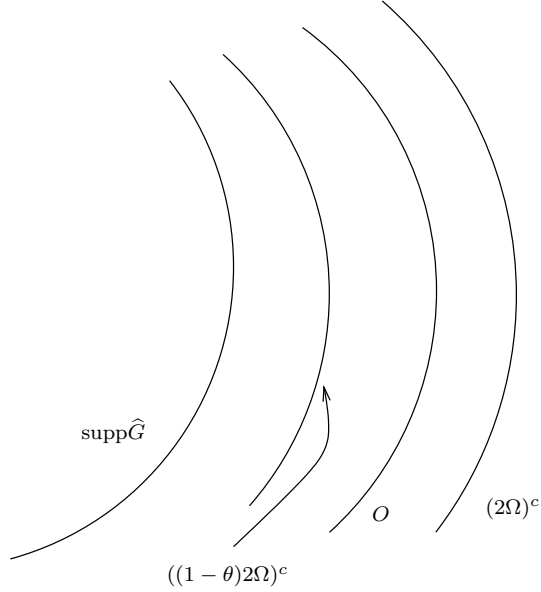}\end{center}
\caption{The domains used in the proof of Theorem \ref{th:main}}
\end{figure}

Here we aim to apply Lemma \ref{th:disjoint-supports},
with $O = \left(\overline{(2-\theta)\Omega}\right)^c$. First let us observe
the following simple fact.

\begin{lemma}
Let $0\le\alpha\le\beta$. 
Then for any bounded convex domain $\Omega$ with $0\in\Omega$ we have 
${\rm dist}\{\alpha\Omega,(\beta\Omega)^c\}=r(\beta-\alpha)$, where 
$r:={\rm dist}\{0,\Omega^c\}=\max \{\rho : B_{\rho} \subseteq \Omega\}$
is the inradius of $\Omega$.
\end{lemma}

Thus we have $\dist\{O,(1-\theta)2\Omega\}=\theta r$ by the above Lemma, hence
for any $\delta<\theta r$ we have $O+B_{\delta}\subset ((1-\theta)2\Omega)^c$.
On the other hand, $\dist\{(2-\theta)\Omega,(2\Omega)^c\}=\theta r >0$ 
again by
the above Lemma and thus $(2\Omega)^c\subset
\left(\overline{(2-\theta)\Omega}\right)^c=O$. 

That is, by (4), 
$\supp \ft{\delta_{\Lambda}}\setminus \{0\} \subseteq (2\Omega)^c \subset O$. 
With these the condition (6) is seen to be fulfilled by $g=G$, as
$\{\ft{g}=0\}=\{\ft{G}=0\} \supseteq ((1-\theta)2\Omega)^c$
since $\supp \ft{G} \subseteq (1-\theta)2\Omega$.

Note that (5) is also satisfied here as 
$\ft{g}(0) = \ft{G}(0) = F(0) =1 $ 
by the definition of $\ft{G}$.
Thus Theorem 5 can be applied and we find that $G$ 
tiles $\RR^d$ with translation set $\Lambda$ at level 
$\ft{g}(0) \ft{\delta_{\Lambda}}\{(0)\} = \ft{G}(0) \dens \Lambda$
by Theorem 6. Here $\ft{G}(0)=1$ and also $\dens \Lambda =1$ in view 
of the considerations next to formula (2). Thus the level of tiling by 
$G$ and $\Lambda$ is found to be 1.
However,
$G(0)=\int \ft{G} > 1$, hence the continuous nonnegative function $G$
can not tile $\RR^d$ at level 1. This contradiction proves that 
there is no function $F$ with the given properties as supposed at the outset.
That is, $\Omega$ is a Tur\'an domain.
\Qed

%%%%%%%%%%%%%%%%%%%%%%%%%%%%%%%%%%%%%%%%%%%%%
% NEW PROOF FOR THE BALL
%%%%%%%%%%%%%%%%%%%%%%%%%%%%%%%%%%%%%%%%%%%%%
\mysection{Another proof that the ball is a Tur\' an domain.}

Here we give a new proof, rather different from that in \cite{ball},
that the ball is a Tur\'an domain. Let $B$ denote the unit ball in $\RR^d$.
To say that the ball is Tur\' an is to prove the inequality
\beql{turan-ineq}
\int f \le 2^{-d} \Abs{B} f(0)
\eeq
for every positive definite $f$ supported in $B$.
By an easy approximation argument it is enough to prove \eqref{turan-ineq}
under the extra assumption that $f$ is smooth.
Noticing further that both sides of the inequality are linear functionals
of $f$ invariant under rotation, we can assume that $f$ is radial
by examining the spherical average of $f$
$$
\int f(Tx)~dT
$$
the integral being over $T \in O_n$, the group of all orthogonal
transformations equipped with Haar measure.

The key ingredient in the proof is the following result.
\begin{othm}\label{th:rudin}
{\rm ({\bf Rudin \cite{rudin}})}\\
Suppose $f$ is a radial smooth positive definite function with
support in the ball $B$. Then $f$ can be written as
a uniformly convergent series
\beql{series}
f = \sum_{k=1}^\infty f_k*\widetilde{f_k},
\ \ \ (\widetilde{f_k}(x) = \overline{f_k(-x)}),
\eeq
with $f_k$ being smooth and supported in the half ball $(1/2)B$.
\end{othm}

Notice that for any integrable function $g$ with support in the
compact set $K$ we have for $f = g*\widetilde{g}$
\begin{eqnarray*}
\int f &=& \Abs{\int g}^2\\
	&\le & \int\Abs{g}^2\Abs{K}\ \ \ \ \mbox{(Cauchy-Schwartz)}\\
	&\le & \int\Abs{g}^2 \cdot 2^{-d} \Abs{K-K}\ \ \ \ \mbox{(Brunn-Minkowski)}\\
	&=& f(0) \cdot 2^{-d} \Abs{K-K}. 
\end{eqnarray*}

Taking $K = \overline{(1/2)B}$ we obtain for every $f_k$ in
\eqref{series}:
$$
\int f_k*\widetilde{f_k} \le 2^{-d}\Abs{B} (f_k*\widetilde{f_k})(0).
$$
Summing the series we have \eqref{turan-ineq} for $f$.
\Qed

%%%%%%%%%%%%%%%%%%%%%%%%%%%%%%%%%%%%%%%%%%%%%
% REFERENCES
%%%%%%%%%%%%%%%%%%%%%%%%%%%%%%%%%%%%%%%%%%%%%
\noindent
\ \\
{\bf Bibliography}
\\

\noindent
{\sc\small
Department of Mathematics, University of Crete, Knossos Ave., \\
714 09 Iraklio, Greece.\\
E-mail: {\tt mk@fourier.math.uoc.gr}\\
\ \\
and\\
\ \\
{\sc\small
Alfr\' ed R\' enyi Institute of Mathematics, Hungarian Academy of Sciences, \\
1364 Budapest, Hungary}\\
E-mail: {\tt revesz@renyi.hu}
}

\end{document}